# Determinantal Processes and Independence[*]

## J. Ben Hough[†]

*Department of Mathematics, U.C. Berkeley, CA 94720, USA*
*e-mail:* `jbhough@math.berkeley.edu`
*url:* `www.math.berkeley.edu/~jbhough`

## Manjunath Krishnapur[†]

*Department of Statistics, U.C. Berkeley, CA 94720, USA*
*e-mail:* `manju@stat.berkeley.edu`
*url:* `stat-www.berkeley.edu/~manju`

## Yuval Peres[†]

*Departments of Statistics and Mathematics, U.C. Berkeley, CA 94720, USA*
*e-mail:* `peres@stat.berkeley.edu`
*url:* `stat-www.berkeley.edu/~peres`

## and

## Bálint Virág[‡]

*Departments of Mathematics and Statistics, University of Toronto, ON, M5S 3G3, Canada*
*e-mail:* `balint@math.toronto.edu`
*url:* `www.math.toronto.edu/~balint`

**Abstract:** We give a probabilistic introduction to determinantal and permanental point processes. Determinantal processes arise in physics (fermions, eigenvalues of random matrices) and in combinatorics (nonintersecting paths, random spanning trees). They have the striking property that the number of points in a region $D$ is a sum of independent Bernoulli random variables, with parameters which are eigenvalues of the relevant operator on $L^2(D)$. Moreover, any determinantal process can be represented as a mixture of determinantal projection processes. We give a simple explanation for these known facts, and establish analogous representations for permanental processes, with geometric variables replacing the Bernoulli variables. These representations lead to simple proofs of existence criteria and central limit theorems, and unify known results on the distribution of absolute values in certain processes with radially symmetric distributions.

Received May 2005.

---

[*]This is an original survey paper.
[†]Research supported in part by NSF grants #DMS-0104073 and #DMS-0244479.
[‡]Research supported in part by grants from the Canada Research Chairs program, NSERC and by the MSRI program on Probability, Algorithms and Statistical Physics.

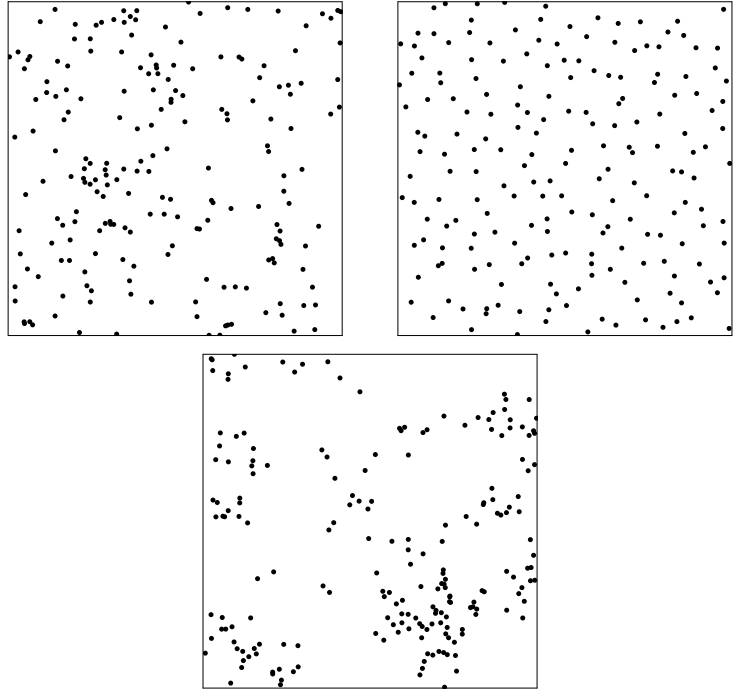

FIG 1. *Samples of translation invariant point processes in the plane: Poisson (left), determinantal (center) and permanental for $K(z,w) = \frac{1}{\pi}e^{z\overline{w}-\frac{1}{2}(|z|^2+|w|^2)}$. Determinantal processes exhibit repulsion, while permanental processes exhibit clumping.*

## 1. Introduction

Determinantal point processes were first studied in 1975 by Macchi [22], who was motivated by fermions in quantum mechanics. Determinantal processes arise naturally in several other settings, including eigenvalues of random matrices, random spanning trees and nonintersecting paths; see, e.g., Burton and Pemantle [4], Soshnikov [27], Lyons [20], Lyons and Steif [21], Shirai and Takahashi [25], Johansson [14], Borodin, Okounkov and Olshanski [3], and Diaconis [8]. A determinantal point process, on a Polish space $\Lambda$ (assumed locally compact) with a reference measure (assumed Radon) $\mu$, is determined by a kernel $K(x,y)$: the joint intensities of the process can be written as $\det(K(x_i,x_j))$. The kernel defines an integral operator $\mathcal{K}$ acting on $L^2(\Lambda)$ that is assumed to be self-adjoint, non-negative and locally trace class, i.e. for every compact $D$ the eigenvalues, $\{\lambda_i^D\}$, of the operator $\mathcal{K}$ restricted to $D$ satisfy $\sum_i \lambda_i^D < \infty$. Determinantal point processes have a special property (Shirai and Takahashi [26] Proposition 2.8) that seems to have only been used in special cases ([1], [23]):

> *In a determinantal process, the number of points that fall in a compact set $D \subset \Lambda$, has the same distribution as a sum of independent Bernoulli($\lambda_i^D$) random variables.*

The proof is immediate from well known formulas for the generating function of particle counts. However, we give a proof starting from first principles and avoid the use of Fredholm determinants.

Permanental processes are defined analogously and are the counterpart of determinantal processes for modeling bosons (see Macchi [22]). They fall into the more general class known as Cox processes [6]. In this case we have:

*In a permanental process, the number of points that fall in a compact set $D \subset \Lambda$, has the same distribution as a sum of independent geometric($\frac{\lambda_i^D}{\lambda_i^D+1}$) random variables.*

For the examples of interest to us, operator theoretic nuances are not essential. The reader will not miss anything significant by keeping in mind just the following cases.

1. $\Lambda$ is a finite set and $K$ is a Hermitian non-negative definite $|\Lambda| \times |\Lambda|$ matrix, and $\mu$ is the counting measure on $\Lambda$.
2. $\Lambda$ is an open set in $R^d$, $\mu$ is Lebesgue measure and $K(x,y)$ is continuous function defining a self-adjoint non-negative integral operator $\mathcal{K}$ on $L^2(\Lambda)$.

A **point process** in a locally compact Polish space $\Lambda$ is a random integer-valued positive Radon measure $\mathcal{X}$ on $\Lambda$. (Recall that a Radon measure is a Borel measure which is finite on compact sets.) If $\mathcal{X}$ almost surely assigns at most measure 1 to singletons, it is a **simple** point process; in this case $\mathcal{X}$ can be identified with a random discrete subset of $\Lambda$, and $\mathcal{X}(D)$ represents the number of points of this set that fall in $D$.

The distribution of a point process can, in most cases, be described by its **joint intensities** (also known as correlation functions).

**Definition 1.** The joint intensities of a point process $\mathcal{X}$ w.r.t. $\mu$ are functions (if any exist) $\rho_k : \Lambda^k \to [0, \infty)$ for $k \geq 1$, such that for any family of mutually disjoint subsets $D_1, \ldots, D_k$ of $\Lambda$,

$$\mathbf{E}\left[\prod_{i=1}^{k} \mathcal{X}(D_i)\right] = \int_{\prod_i D_i} \rho_k(x_1, \ldots, x_k) d\mu(x_1) \ldots d\mu(x_k), \qquad (1)$$

**Remark 2.** For overlapping sets, the situation is more complicated. Restricting attention to simple point processes, $\rho_k$ is not the intensity measure of $\mathcal{X}^k$, but that of $\mathcal{X}^{\wedge k}$, the set of ordered $k$-tuples of distinct points of $\mathcal{X}$. Indeed, (1) implies (see [18, 19, 23]) that for any Borel set $B \subset \Lambda^k$ we have

$$\mathbf{E}\,\#(B \cap \mathcal{X}^{\wedge k}) = \int_B \rho_k(x_1, \ldots, x_k)\, d\mu(x_1) \ldots d\mu(x_k)\,. \qquad (2)$$

When $B = \prod D_i^{\otimes k_i}$ for a mutually disjoint family of subsets $D_1, \ldots, D_r$ of $\Lambda$, and $k = \sum_{i=1}^{r} k_i$, the left hand side becomes

$$\mathbf{E}\left[\prod_{i=1}^{r} \binom{\mathcal{X}(D_i)}{k_i} k_i!\right]. \qquad (3)$$

For a general point process $\mathcal{X}$, observe that it can be identified with a simple point process $\mathcal{X}^*$ on $\Lambda \times \{1, 2, 3, \ldots\}$ such that $\mathcal{X}^*(D \times \{1, 2, 3, \ldots\}) = \mathcal{X}(D)$ for Borel $D \subset \Lambda$. Thus, if $\mathcal{X}(D)$ has exponential tails for all compact $D \subset \Lambda$, then the joint intensities determine the law of $\mathcal{X}$; see [18, 19]. In particular, Theorems 7 and 10 below imply that this is the case for determinantal and permanental processes governed by a trace class operator, since a convergent sum of Bernoulli (or geometric) variables always has exponential tails.

**Definition 3.** A point process $\mathcal{X}$ on $\Lambda$ is said to be a **determinantal process** with kernel $K$ if it is simple and its joint intensities satisfy:

$$\rho_k(x_1, \ldots, x_k) = \det\left(K(x_i, x_j)\right)_{1 \le i, j \le k}, \tag{4}$$

for every $k \geq 1$ and $x_1, \ldots, x_k \in \Lambda$.

**Remark 4.** We postulate determinantal processes to be simple because we have adopted equation (1) as the definition of joint intensities. If instead, we start with the slightly more restrictive definition of joint intensities as explained above in Remark 2, it follows that the process must be simple. Nevertheless, postulating simplicity is more in tune with the conventions of Physicists who often consider a determinantal process with $k$ points, not as a random counting measure but as a random point in $\Lambda^k/\text{Diagonals}$, where "Diagonals" denotes the subset of points of $\Lambda^k$ with at least two co-ordinates equal. This viewpoint together with a postulate on the behaviour of quantum amplitudes under exchange of particles leads naturally to Determinantal and Permanental processes (and additionally to "fractional statistics", when $\Lambda$ is two dimensional, see for instance [16]).

Consider a kernel $K$ that defines a self-adjoint trace-class operator $\mathcal{K}$. Macchi [22] and Soshnikov [27] showed that there exists a determinantal point process with kernel $K$ if and only if all eigenvalues of $\mathcal{K}$ are in the interval $[0, 1]$. In Section 2 we shall give a probabilistic proof of this fact.

**Remark 5.** When we speak of a kernel $K$ on $\Lambda^2$, a priori it is only defined only almost everywhere w.r.t. $\mu \times \mu$ and thus quantities like $\int K(x, x) d\mu(x)$ that appear in the definition of joint intensities are not defined. This can be made sense of as follows.

The kernel $K$ defines a self-adjoint integral operator $\mathcal{K}$ that has eigenfunctions $\phi_k$ and eigenvalues $\lambda_k$, where $\phi_k$ are orthogonal in $L^2(\mu)$. In particular, there is a set $\Lambda_1 \subset \Lambda$ such that $\phi_k$ are all defined point-wise on $\Lambda_1$ and $\mu(\Lambda_1^c) = 0$. At least in the case when $\mathcal{K}$ has finite rank, this shows that $K(x, y) = \sum_k \lambda_k \phi_k(x) \overline{\phi}_k(y)$ is well defined on $\Lambda_1 \times \Lambda_1$, and that is sufficient to define $\int K(x, x) d\mu(x)$ etc. For more details on this point, see Lemmas 1,2 of the survey paper of Soshnikov [27].

Recall that the permanent of an $n \times n$ matrix $M$ is defined as

$$\text{per}(M) = \sum_{\pi \in S_n} \prod_{i=1}^n M_{i, \pi(i)}. \tag{5}$$

**Definition 6.** A point process $\mathcal{X}$ on $\Lambda$ is said to be a **permanental process** with kernel $K$ if its joint intensities satisfy:

$$\rho_k(x_1,\ldots,x_k) = \operatorname{per}\left(K(x_i,x_j)\right)_{1\leq i,j\leq k}, \qquad (6)$$

for every $k \geq 1$ and $x_1,\ldots,x_k \in \Lambda$.

For any kernel $K$ that defines a self-adjoint non-negative definite operator $\mathcal{K}$ (i.e., the eigenvalues of $\mathcal{K}$ are nonnegative), there exists a permanental point process with kernel $K$. We shall give a proof of this known fact in Section 4.

Now we state the main theorems. The most common application of the following theorem is to describe the behavior of a determinantal process already restricted to a subset.

**Theorem 7.** *Suppose $\mathcal{X}$ is a determinantal process with a trace-class kernel $K$. Write*

$$K(x,y) = \sum_{k=1}^{n} \lambda_k \phi_k(x) \overline{\phi}_k(y),$$

*where $\phi_k$ are normalized eigenfunctions of $K$ with eigenvalues $\lambda_k \in [0,1]$. (Here $n = \infty$ is allowed). Let $I_k$, $1 \leq k \leq n$ be independent random variables with $I_k \sim \operatorname{Bernoulli}(\lambda_k)$. Set*

$$K_I(x,y) = \sum_{k=1}^{n} I_k \phi_k(x) \overline{\phi}_k(y).$$

*$K_I$ is a random analogue of the kernel $K$. Let $\mathcal{X}_I$ be the determinantal process with kernel $K_I$ (i.e., first choose the $I_k$'s and then independently sample a discrete set that is determinantal with kernel $K_I$). Then*

$$\mathcal{X} \stackrel{d}{=} \mathcal{X}_I. \qquad (7)$$

*In particular, the total number of points in the process $\mathcal{X}$ has the distribution of a sum of independent $\operatorname{Bernoulli}(\lambda_k)$ random variables.*

In the special case of random spanning trees of a finite graph, Bapat [1] was the first to observe the last fact stated above. Namely, he proved that the number of edges of the spanning tree falling in a subset of edges of the given graph has the distribution of a sum of independent Bernoullis.

When the $D_i$'s are related in a special way, there exists a simple probabilistic description of the joint distributions of the counts $\mathcal{X}(D_i)$.

**Definition 8.** Given an integral kernel $K$ acting on $L^2(\Lambda)$, the subsets $D_1,\ldots,D_r$ of $\Lambda$ with $D = \cup_i D_i$ are **simultaneously observable** if the eigenfunctions of the restricted kernel $K|_{D\times D}$ are also eigenfunctions of $K|_{D_i \times D_i}$ acting on $L^2(D_i)$ for every $1 \leq i \leq r$.

The motivation for this terminology comes from quantum mechanics, where two physical quantities can be simultaneously measured if the corresponding operators commute.

**Proposition 9.** *Under the assumptions of Theorem 7, let $D_i \subset \Lambda$, $1 \leq i \leq r$ be mutually disjoint and simultaneously observable. Let $\mathbf{e}_i$ be the standard basis vectors in $\mathbb{R}^r$. Denote by $\phi_k$, the common eigenfunctions of $K$ on the $D_i$'s and by $\lambda_{k,i}$ the corresponding eigenvalues. Write $\lambda_k = \sum_i \lambda_{k,i}$ and note that $\lambda_k \leq 1$. Then*

$$(\mathcal{X}(D_1),\ldots,\mathcal{X}(D_r)) \stackrel{d}{=} \sum_k (\xi_{k,1},\ldots,\xi_{k,r}), \qquad (8)$$

*where $\vec{\xi}_k = (\xi_{k,1},\ldots,\xi_{k,r})$ are independent for different values of $k$, with $\mathbf{P}(\vec{\xi}_k = \mathbf{e}_i) = \lambda_{k,i}$ for $1 \leq i \leq r$ and $\mathbf{P}(\vec{\xi}_k = \mathbf{0}) = 1 - \lambda_k$. In words, $(\mathcal{X}(D_1),\ldots,\mathcal{X}(D_r))$ has the same distribution as the vector of counts in $r$ cells, if we pick $n$ balls and assign the $k^{th}$ ball to the $i^{th}$ cell with probability $\lambda_{k,i}$ (there may be a positive probability of not assigning it to any of the cells).*

**Theorem 10.** *Suppose $\mathcal{X}$ is a permanental process in $\Lambda$ with a trace-class kernel $K$. Write $K(x,y) = \sum_{k=1}^n \lambda_k \phi_k(x) \overline{\phi}_k(y)$, where $\phi_k$ are normalized eigenfunctions of $K$ with eigenvalues $\lambda_k$ ($n = \infty$ is allowed). Let $\vec{\alpha} = (\alpha_1,\ldots,\alpha_n)$, where $\alpha_i$ are non-negative integers such that $\ell = \ell(\vec{\alpha}) = \alpha_1 + \cdots + \alpha_n < \infty$ and let $Z^{\vec{\alpha}}$ be the random vector in $\mathbb{R}^\ell$ with density:*

$$p(z_1,\ldots,z_\ell) = \frac{1}{\ell! \alpha_1! \cdots \alpha_n!} \left| \operatorname{per} \begin{bmatrix} \{\phi_1(z_1) & \cdots & \phi_1(z_\ell)\} \alpha_1 \\ \cdots & \cdots & \cdots \\ \{\phi_n(z_1) & \cdots & \phi_n(z_\ell)\} \alpha_n \end{bmatrix} \right|^2, \qquad (9)$$

*where the notation $\{\phi_i(z_1) \cdots \phi_i(z_\ell)\} \alpha_i$ indicates that the row $\phi_i(z_1) \cdots \phi_i(z_\ell)$ is repeated $\alpha_i$ times. Let $\mathcal{X}_{\vec{\alpha}}$ be the point process determined by $Z^{\vec{\alpha}}$, i.e., $\mathcal{X}_{\vec{\alpha}}(D)$ is the number of $j \leq \ell$ such that $Z_j^{\vec{\alpha}} \in D$. Let $\gamma_1,\ldots,\gamma_n$ be independent geometric random variables with $\mathcal{P}(\gamma_k = s) = \left(\frac{\lambda_k}{\lambda_k + 1}\right)^s \left(\frac{1}{\lambda_k + 1}\right)$, for $s = 0,1,2,\ldots$. Then*

$$\mathcal{X} \stackrel{d}{=} \mathcal{X}_{\vec{\gamma}},$$

*where $\vec{\gamma} = (\gamma_1,\ldots,\gamma_n)$. In particular, $\mathcal{X}(\Lambda)$ has the distribution of a sum of independent geometric$(\frac{\lambda_k}{\lambda_k+1})$ random variables.*

**Remark 11.** The density given in (9) has physical significance. Interpreting the functions $\phi_k$ as eigenstates of a one-particle Hamiltonian, (9) gives the distribution for $\ell$ non-interacting bosons in a common potential given that $\alpha_i$ of them lie in the eigenstate $\phi_i$. This density is the exact analogue of the density

$$p(z_1,\ldots,z_\ell) = \frac{1}{\ell!} \left| \det \begin{bmatrix} \phi_{i_1}(z_1) & \cdots & \phi_{i_1}(z_\ell) \\ \cdots & \cdots & \cdots \\ \phi_{i_n}(z_1) & \cdots & \phi_{i_n}(z_\ell) \end{bmatrix} \right|^2, \qquad (10)$$

which appears in Theorem 7 and gives the distribution for $\ell$ non-interacting fermions in a common potential given that one fermion lies in each of the eigenstates $\phi_{i_1},\ldots,\phi_{i_\ell}$. The fact that (10) vanishes if a row is repeated illustrates Pauli's exclusion principle, which states that multiple fermions cannot occupy the same eigenstate. See [9] for more details.

**Theorem 12.** *Under the assumptions of Theorem 10, suppose $D_1, \ldots, D_r$ are simultaneously observable as in definition 8. Denote by $\phi_k$ the common eigenfunctions of $K$ on the $D_i$'s and by $\lambda_{k,i}$ the corresponding eigenvalues. Then*

$$(\mathcal{X}(D_1), \ldots, \mathcal{X}(D_r)) \stackrel{d}{=} \sum_k (\eta_{k,1}, \ldots, \eta_{k,r}), \tag{11}$$

*where $(\eta_{k,1}, \ldots, \eta_{k,r})$ are independent for different values of $k$, for each $k$, the sum $\eta_k = \sum_i \eta_{k,i}$ has a geometric distribution with mean $\lambda_k := \sum_i \lambda_{k,i}$ and given $\sum_i \eta_{k,i} = N$,*

$$(\eta_{k,1}, \ldots, \eta_{k,r}) \stackrel{d}{=} \text{Multinomial}\left(N; \frac{\lambda_{k,1}}{\lambda_k}, \ldots, \frac{\lambda_{k,r}}{\lambda_k}\right).$$

## 2. Determinantal processes

We begin with a few important examples of determinantal processes.

**Example 13 (Non-intersecting random walks).** Consider $n$ independent simple symmetric random walks on $\mathbb{Z}$ started from $i_1 < i_2 < \ldots < i_n$ where all the $i_j$'s are even. Let $P_{i,j}(t)$ be the $t$-step transition probabilities. Karlin and McGregor [15] show that the probability that the random walks are at $j_1 < j_2 < \ldots < j_n$ at time $t$ and have mutually disjoint paths is

$$\det \begin{pmatrix} P_{i_1,j_1}(t) & \ldots & P_{i_1,j_n}(t) \\ \ldots & \ldots & \ldots \\ P_{i_n,j_1}(t) & \ldots & P_{i_n,j_n}(t) \end{pmatrix}.$$

It follows easily that if $t$ is even and we also condition the walks to return to $i_1, \ldots, i_n$ at time $t$, then the positions of the walks at time $t/2$ are determinantal with a Hermitian kernel. See Johansson [14] for this and more general results.

**Example 14 (Uniform spanning trees).** Let $G$ be a finite undirected graph and let $E$ be the set of oriented edges (each undirected edge appears in $E$ with both orientations). Let $T$ be uniformly chosen from the set of spanning trees of $G$. For each directed edge $e = vw$, let $\chi^e := \mathbf{1}_{vw} - \mathbf{1}_{wv}$ denote the unit flow along $e$. Define

$$\begin{aligned}
\mathbf{H} &= \{f : E \to \mathbb{R} : f(vw) = -f(wv) \, \forall v, w\} \\
\bigstar &= \text{span}\{\sum_w \chi^{vw} : \text{ where } v \text{ is a vertex.}\} \\
\diamondsuit &= \text{span}\{\sum_{i=1}^n \chi^{e_i} : e_1, \ldots, e_n \text{ is an oriented cycle}\}.
\end{aligned}$$

It is easy to see that $\mathbf{H} = \bigstar \oplus \diamondsuit$. Now, define $I^e := \mathcal{K}_\bigstar \chi^e$, the orthogonal projection onto $\bigstar$. Then, the set of edges in $T$ forms a determinantal process with kernel $K(e, f) := (I^e, I^f)$ with respect to counting measure on the set of edges. This was proved by Burton and Pemantle [4], who represented $K(e, f)$ as the current flowing through $f$ when a unit of current is sent from the tail to the head of $e$. The Hilbert space formulation above is from BLPS [2].

**Example 15 (Ginibre ensemble).** Let $Q$ be an $n \times n$ matrix with i.i.d. standard complex normal entries. Ginibre [10] proved that the eigenvalues of $Q$ form a determinantal process in $\mathbb{C}$ with the kernel

$$K_n(z, w) = \frac{1}{\pi} e^{-\frac{1}{2}(|z|^2+|w|^2)} \sum_{k=0}^{n-1} \frac{(z\overline{w})^k}{k!},$$

with respect to Lebesgue measure. As $n \to \infty$, we get a determinantal process in the plane with the kernel

$$\begin{aligned} K(z, w) &= \frac{1}{\pi} e^{-\frac{1}{2}(|z|^2+|w|^2)} \sum_{k=0}^{\infty} \frac{(z\overline{w})^k}{k!} \\ &= \frac{1}{\pi} e^{-\frac{1}{2}(|z|^2+|w|^2)+z\overline{w}}. \end{aligned}$$

**Example 16 (Zero set of a Gaussian analytic function).** The power series $\mathbf{f}_1(z) := \sum_{n=0}^{\infty} a_n z^n$ where $a_n$ are i.i.d. standard complex normals, defines a random analytic function in the unit disk (almost surely). Peres and Virág [23] show that the zero set of $\mathbf{f}_1$ is a determinantal process in the disk with the Bergman kernel

$$K(z, w) = \frac{1}{\pi(1 - z\overline{w})^2} = \frac{1}{\pi} \sum_{k=0}^{\infty} (k+1)(z\overline{w})^k,$$

with respect to Lebesgue measure in the unit disk.

### *Determinantal projection processes: motivation and construction*

The most general determinantal processes are mixtures of *determinantal projection processes*, i.e. processes whose kernel $K_H$ defines a projection operator $\mathcal{K}_H$ to a subspace $H \subset L^2(\Lambda)$ or, equivalently, $K_H(x, y) = \sum_{k=1}^{n} \phi_k(x)\overline{\phi}_k(y)$ where $\{\phi_k\}$ is any orthonormal basis for $H$.

**Lemma 17.** *Suppose $\mathcal{X}$ is a determinantal projection process on $\Lambda$, with kernel $K(x, y) = \sum_{k=1}^{n} \phi_k(x)\overline{\phi}_k(y)$ where $\{\phi_k\}$ is an orthonormal set in $L^2(\Lambda)$. Then the number of points in $\mathcal{X}$ is equal to $n$, almost surely.*

*Proof.* The conditions imply that the matrix $(K(x_i, x_j))_{1 \leq i,j \leq k}$ has rank at most $n$ for any $k \geq 1$. From (3), we see that $\mathbf{E}\left[\binom{\mathcal{X}(\Lambda)}{k}\right] = 0$ for $k > n$. This shows that $\mathcal{X}(\Lambda) \leq n$ almost surely. However,

$$\begin{aligned}
\mathbf{E}\left[\mathcal{X}(\Lambda)\right] &= \int_\Lambda \rho_1(x)\,d\mu(x) \\
&= \int_\Lambda K(x,x)\,d\mu(x) \\
&= \sum_{k=1}^n \int_\Lambda |\phi_k(x)|^2 d\mu(x) \\
&= n.
\end{aligned}$$

Therefore $\mathcal{X}(\Lambda) = n$, almost surely. $\square$

Despite the fact that determinantal processes arise naturally and many important statistics can be computed, the standard Definition 3 is lacking in direct probabilistic intuition. Below we present an algorithm that is somewhat more natural from a probabilist's point of view, and can also be used for modeling determinantal processes.

In the discrete case, the projection operator $\mathcal{K}_H$ can be applied to the delta function at a point, and we have $\mathcal{K}_H \delta_x(\cdot) = K(\cdot, x)$. In the general case we take this as the *definition* of $\mathcal{K}_H \delta_x$. Let $\|\cdot\|$ denote the norm of $L^2(\mu)$. The intensity measure of the process is given by

$$\mu_H(x) = \rho_1(x) d\mu(x) = \|\mathcal{K}_H \delta_x\|^2 d\mu(x). \tag{12}$$

When $\mu$ is supported on countably many points, we have $\|\mathcal{K}_H \delta_x\| = \mathrm{dist}(\delta_x, H^\perp)$ (where $\perp$ denotes orthocomplement), giving a natural interpretation of the intensity $\rho_1$ in general.

Note that $\mu_H(\Lambda) = \dim(H)$, so $\mu_H/\dim(H)$ is a probability measure on $\Lambda$. We construct the determinantal process as follows. Start with $n = \dim(H)$, and $H_n = H$.

**Algorithm 18.**
- If $n = 0$, stop.
- Pick a random point $X_n$ from the probability measure $\mu_{H_n}/n$.
- Let $H_{n-1} \subset H_n$ be the orthocomplement of the function $\mathcal{K}_{H_n}\delta_x$ in $H_n$. In the discrete case, $H_{n-1} = \{f \in H_n : f(X_n) = 0\}$. Note that $\dim(H_{n-1}) = n - 1$ a.s.
- Decrease $n$ by 1 and iterate.

**Proposition 19.** *The points $(X_1, \ldots, X_n)$ constructed by Algorithm 18 are distributed as a uniform random ordering of the points in a determinantal process $\mathcal{X}$ with kernel $K$.*

*Proof.* Let $\psi_j = \mathcal{K}_H \delta_{x_j}$. Projecting to $H_j$ is equivalent to first projecting to $H$ and then to $H_j$, and it is easy to check that $\mathcal{K}_{H_j} \delta_{x_j} = \mathcal{K}_{H_j} \psi_j$. Thus, by (12), the density of the random vector $(X_1, \ldots, X_n)$ constructed by the algorithm equals

$$p(x_1, \ldots, x_n) = \prod_{j=1}^n \frac{\|\mathcal{K}_{H_j} \psi_j\|^2}{j}.$$

Note that $H_j = H \cap \langle \psi_{j+1}, \ldots, \psi_n \rangle^\perp$, and therefore $V = \prod_{j=1}^n \|\mathcal{K}_{H_j}\psi_j\|$ is exactly the repeated "base times height" formula for the volume of the parallelepiped determined by the vectors $\psi_1, \ldots, \psi_n$ in the finite-dimensional vector space $H \subset L^2(\Lambda)$. It is well-known that $V^2$ equals the determinant of the *Gram matrix* whose $i,j$ entry is given by the scalar product of $\psi_i, \psi_j$, that is $\int \psi_i \overline{\psi_j} d\mu = K(x_i, x_j)$. We get

$$p(x_1, \ldots, x_n) = \frac{1}{n!} \det(K(x_i, x_j)),$$

so the random variables $X_1, \ldots, X_n$ are exchangeable. Viewed as a point process, the $n$-point joint intensity of $\{X_j\}_{j=1}^n$ is $n! p(x_1, \ldots, x_n)$, which agrees with that of the determinantal process $\mathcal{X}$. The claim now follows by Lemma 17. □

**Example 20 (Uniform spanning trees).** We continue the discussion of Example 14 Let $G_{n+1}$ be an undirected graph on $n+1$ vertices. For every edge $e$, the effective resistance when a unit of current is sent along $e$ is $R(e) = (I^e, I^e)$. To use our algorithm to choose a uniform spanning tree, proceed as follows:

- If $n = 0$, stop.
- Take $X_n$ to be a random edge, chosen so that $\mathbf{P}(X_n = e_i) = \frac{1}{n} R(e_i)$.
- Construct $G_n$ from $G_{n+1}$ by contracting the edge $X_n$, and update the effective resistances $\{R(e)\}$.
- Decrease $n$ by one and iterate.

For sampling uniform spanning trees, more efficient algorithms are known, but for the general case the above procedure is the most efficient we are aware of.

We shall need the following lemmas.

**Lemma 21.** *Suppose $\{\phi_k\}_{k=1}^n$ is an orthonormal set in $L^2(\Lambda)$. Then there exists a determinantal process with kernel $K(x,y) = \sum_{k=1}^n \phi_k(x)\overline{\phi}_k(y)$.*

*Proof.* For any $x_1, \ldots, x_n$ we have $(K(x_i, x_j))_{1 \le i, j \le n} = A A^*$, where $A_{i,k} = \phi_k(x_i)$. Therefore, $\det(K(x_i, x_j))$ is non-negative. Moreover,

$$\int_{\Lambda^n} \det(K(x_i, x_j))_{i,j} dx_1 \ldots dx_n$$

$$= \int_{\Lambda^n} \det(\phi_j(x_i))_{i,j} \det(\overline{\phi}_i(x_j))_{i,j} dx_1 \ldots dx_n$$

$$= \int_{\Lambda^n} \sum_{\pi, \tau \in S_n} \mathrm{sgn}(\pi)\mathrm{sgn}(\tau) \prod_{k=1}^n \phi_{\pi(k)}(x_k)\overline{\phi}_{\tau(k)}(x_k) dx_1 \ldots dx_n.$$

In the sum, if $\pi(k) \ne \tau(k)$, then $\int_\Lambda \phi_{\pi(k)}(x_k)\overline{\phi}_{\tau(k)}(x_k)dx_k = 0$, and when $\pi(k) = \tau(k)$, this integral is 1. Thus, only the terms with $\pi = \tau$ contribute. We get

$$\int_{\Lambda^n} \det(K(x_i, x_j))_{1 \le i,j \le n} dx_1 \ldots dx_n = n!,$$

which along with the non-negativity of $\det\left(K(x_i, x_j)\right)_{1 \leq i,j \leq n}$ shows that $\frac{1}{n!}\det\left(K(x_i, x_j)\right)_{1 \leq i,j \leq n}$ is a probability density on $\Lambda^n$. If we look at the resulting random variable as a set of unlabeled points in $\Lambda$, we get the desired $n$-point joint intensity $\rho_n$.

Lower joint intensities are got by integrating over some of the $x_i$s:

$$\rho_k(x_1,\ldots,x_k) = \frac{1}{(n-k)!}\int_{\Lambda^{n-k}} \rho_n(x_1,\ldots,x_n)\prod_{j>k} d\mu(x_j). \qquad (13)$$

We caution that (13) is valid only for a point process that has $n$ points almost surely. In general, there is no way to get lower joint intensities from higher ones.

We now show how to get $\rho_{n-1}$. The others can be found exactly the same way, or inductively. Set $k = n-1$ in (13) and expand $\rho_n(x_1,\ldots,x_n) = \det\left(K(x_i,x_j)\right)_{1 \leq i,j \leq n}$ as we did before to get

$$\rho_{n-1}(x_1,\ldots,x_{n-1})$$
$$= \sum_{\pi,\tau} \mathrm{sgn}(\pi)\mathrm{sgn}(\tau)\prod_{k=1}^{n-1}\phi_{\pi(k)}(x_k)\overline{\phi}_{\tau(k)}(x_k)\int_\Lambda \phi_{\pi(n)}(x_n)\overline{\phi}_{\tau(n)}(x_n)d\mu(x_n).$$

If $\pi(n) \neq \tau(n)$, the integral vanishes. And if $\pi(n) = \tau(n) = j$, $\pi$ and $\tau$ map $\{1,\ldots,n-1\}$ to $\{1,2,\ldots,n\} - \{j\}$ (and the product of the signs of these "permutations" is the same as $\mathrm{sgn}(\pi)\mathrm{sgn}(\tau)$, because $\pi(n) = \tau(n)$). This gives us

$$\rho_{n-1}(x_1,\ldots,x_{n-1}) = \sum_{j=1}^n \det\left(\phi_k(x_i)\right)_{1 \leq i \leq n-1, k \neq j} \det\left(\overline{\phi}_k(x_i)\right)_{k \neq j, 1 \leq i \leq n-1}.$$

We must show that this quantity is equal to $\det\left(K(x_i,x_j)\right)_{i,j \leq n-1}$. For this note that

$$\left(K(x_i,x_j)\right)_{i,j \leq n-1} = \left(\phi_k(x_i)\right)_{1 \leq i \leq n-1, k \leq n}\left(\overline{\phi}_k(x_i)\right)_{k \leq n, i \leq n-1},$$

and apply the Cauchy-Binet formula. Recall that for matrices $A, B$ of orders $m \times n$ and $n \times m$ respectively, the Cauchy-Binet formula says

$$\det(AB) = \sum_{1 \leq i_1,\ldots,i_m \leq n} \det\left(A[i_1,\ldots,i_m]\right)\det\left(B[i_1,\ldots,i_m]\right), \qquad (14)$$

where we abuse notation and let $A[i_1,\ldots,i_m]$ stand for the matrix formed by taking the *columns* numbered $i_1,\ldots,i_m$ and $B[i_1,\ldots i_m]$ for the matrix formed by the corresponding *rows* of $B$. The identity (14) can be proved by observing that both sides are multi-linear in the rows of $A$ and in the columns of $B$. □

We now prove Theorem 7. Before that we remark that in many examples the kernel $K$ defines a projection operator, i.e, $\lambda_k = 1$ for all $k$. Then $I_k = 1$ for all $k$, almost surely, and the theorem is trivial. The theorem is applicable to the restriction of the process $\mathcal{X}$ to $D$ for any Borel set $D \subset \Lambda$, as the restricted process is determinantal with kernel the restriction of $K$ to $D \times D$.

*Proof of Theorem 7.* First assume that $\mathcal{K}$ is a finite-dimensional operator-i.e., $K(x,y) = \sum_{k=1}^{n} \lambda_k \phi_k(x) \overline{\phi}_k(y)$ for some finite $n$. We show that the processes on the left and right side of (7) have the same joint intensities. By (3), this implies that these processes have the same distribution.

Note that the process $\mathcal{X}_I$ exists by Lemma 21. For $m > n$, the $m$-point joint intensities of both $\mathcal{X}$ and $\mathcal{X}_I$ are clearly zero. Now consider $m \leq n$ and $x_1, \ldots, x_m \in \Lambda$. We claim that:

$$\mathbf{E}\left[\det\left(K_I(x_i, x_j)\right)_{1 \leq i,j \leq m}\right] = \det\left(K(x_i, x_j)\right)_{1 \leq i,j \leq m}. \tag{15}$$

To prove (15), note that

$$(K_I(x_i, x_j))_{1 \leq i,j \leq m} = A\,B, \tag{16}$$

where $A$ is the $m \times n$ matrix with $A_{i,k} = I_k \phi_k(x_i)$ and $B$ is the $n \times m$ matrix with $B_{k,j} = \overline{\phi}_k(x_j)$.

Apply Cauchy-Binet formula (14) to $A, B$ defined above and take expectations. Observe that $B[i_1, \ldots i_m]$ is nonrandom and

$$\mathbf{E}\left[\det\left(A[i_1, \ldots, i_m]\right)\right] = \det\left(C[i_1, \ldots, i_m]\right)$$

where $C$ is the $m \times n$ matrix $C_{i,k} = \lambda_k \phi_k(x_i)$. Applying the Cauchy-Binet formula in the reverse direction to $C$ and $B$, we obtain (15) and hence also (7). By Lemma 17, given $\{I_k\}_{k \geq 1}$, $\mathcal{X}_I$ has $\sum_k I_k$ points, almost surely. Therefore,

$$\mathcal{X}(\Lambda) \stackrel{d}{=} \sum_k I_k.$$

So far we assumed that the operator $\mathcal{K}$ determined by the kernel $K$ is finite dimensional. Now suppose $\mathcal{K}$ is a general trace class operator. Then $\sum_k \lambda_k < \infty$ and hence, almost surely, $\sum_k I_k < \infty$. Therefore, the process $\mathcal{X}_I$ is well defined and (16) is valid by the same reasoning. Taking expectations and observing that the summands in the Cauchy-Binet formula are non-negative, we obtain

$$\mathbf{E}\left[\det\left(K_I(x_i, x_j)\right)_{1 \leq i,j \leq m}\right] = \sum_{1 \leq i_1, \ldots, i_m} \det\left(C[i_1, \ldots, i_m]\right) \det\left(B[i_1, \ldots, i_m]\right),$$

where $C$ is the same as before. To conclude that the right hand side is just $\det\left(K(x_i, x_j)\right)_{1 \leq i,j \leq m}$, we first apply the Cauchy-Binet formula to the finite approximation $(K_N(x_i, x_j))_{1 \leq i,j \leq m}$, where $K_N(x,y) = \sum_{k=1}^{N} \lambda_k \phi_k(x) \overline{\phi}_k(y)$. Letting $N \to \infty$, we see that

$$\mathbf{E}\left[\det\left(K_I(x_i, x_j)\right)_{1 \leq i,j \leq m}\right] = \det\left(K(x_i, x_j)\right)_{1 \leq i,j \leq m},$$

as was required to show. (In short, the proof for the the infinite case is exactly the same as before, only we cautiously avoided applying Cauchy-Binet formula to the product of two infinite rectangular matrices). □

Now we give a probabilistic proof of the following criterion for a Hermitian integral kernel to define a determinantal process.

**Theorem 22 (Macchi [22], Soshnikov [27]).** *Let $K$ determine a self-adjoint integral operator $\mathcal{K}$ on $L^2(\Lambda)$ that is locally trace class. Then $K$ defines a determinantal process on $\Lambda$ if and only if all the eigenvalues of $\mathcal{K}$ are in $[0, 1]$.*

*Proof.* We can assume that $\mathcal{K}$ is trace class, since it suffices to construct a determinantal process on compact subsets of $\Lambda$ with kernel the restriction of $K$.

**Sufficiency:** If $\mathcal{K}$ is a projection operator, this is precisely Lemma 21. If the eigenvalues are $\{\lambda_k\}$, then as in the proof of Theorem 7 we construct the process $\mathcal{X}_I$. The proof there shows that $\mathcal{X}_I$ is determinantal with kernel $K$.

**Necessity:** Suppose that $\mathcal{X}$ is determinantal with kernel $K$. Since the joint intensities of $\mathcal{X}$ are non-negative, $K$ must be non-negative definite. Now suppose that the largest eigenvalue of $\mathcal{K}$ is $\lambda > 1$. Let $\mathcal{X}_1$ be the process obtained by first sampling $\mathcal{X}$ and then independently deleting each point of $\mathcal{X}$ with probability $1 - \frac{1}{\lambda}$. Computing the joint intensities shows that $\mathcal{X}_1$ is determinantal with kernel $\frac{1}{\lambda}K$.

Now $\mathcal{X}$ has finitely many points (we assumed that $\mathcal{K}$ is trace class) and $\lambda > 1$. Hence, $\mathbf{P}[\mathcal{X}_1(\Lambda) = 0] > 0$. However, $\frac{1}{\lambda}K$ has all eigenvalues in $[0, 1]$, with at least one eigenvalue equal to 1, whence by Theorem 7, $P[\mathcal{X}_1(\Lambda) \geq 1] = 1$, a contradiction. □

**Example 23 (Non-measurability of the Bernoullis).** A natural question that arises from Theorem 7 is whether, given a realization of the determinantal process $\mathcal{X}$, we can determine the values of the $I_k$'s. This is not always possible, i.e., the $I_k$'s are not measurable w.r.t. the process $\mathcal{X}$ in general.

Consider the graph $G$ with vertices $\{a, b, c, d\}$ and edges $e_1 = (a, b), e_2 = (b, c), e_3 = (c, d), e_4 = (d, a), e_5 = (a, c)$. By the Burton-Pemantle Theorem [4], the edge-set of a uniformly chosen spanning tree of $G$ is a determinantal process. In this case, the kernel restricted to the set $D = \{e_1, e_2, e_3\}$ is easily computed to be

$$(K(e_i, e_j))_{1 \leq i,j \leq 3} = \begin{pmatrix} 5 & -3 & -1 \\ -3 & 5 & -1 \\ -1 & -1 & -1 \end{pmatrix}.$$

This matrix has eigenvalues $1, \frac{7-\sqrt{17}}{16}, \frac{7+\sqrt{17}}{16}$. Since all measurable events have probabilities that are multiples of $\frac{1}{8}$, it follows that the Bernoullis cannot be measurable.

Theorem 7 gives us the distribution of the number of points $\mathcal{X}(D)$ in any subset of $\Lambda$. Given several regions $D_1, \ldots, D_r$, can we find the joint distribution of $\mathcal{X}(D_1), \ldots, \mathcal{X}(D_r)$? It seems that a simple probabilistic description of the joint distribution exists only when the $D_i$'s are simultaneously observable, as in Theorem 9.

*Proof of Theorem 9.* At first we make the following assumptions:

- $\cup_i D_i = \Lambda$.

- $K$ defines a finite dimensional projection operator-i.e., $K(x, y) = \sum_{k=1}^{n} \phi_k(x)\overline{\phi}_k(y)$ for $x, y \in \Lambda$, and $\{\phi_k\}$ is an orthonormal set in $L^2(\Lambda)$.

Note that by our assumption, $\phi_k$ are also orthogonal on $D_i$ for every $1 \leq i \leq r$. Moreover, it is clear that $\lambda_{k,i} = \int_{D_i} |\phi_k|^2$.

We write

$$\begin{pmatrix} K(x_1, x_1) & \ldots & K(x_1, x_n) \\ \ldots & \ldots & \ldots \\ K(x_n, x_1) & \ldots & K(x_n, x_n) \end{pmatrix}$$
$$= \begin{pmatrix} \phi_1(x_1) & \ldots & \phi_n(x_1) \\ \ldots & \ldots & \ldots \\ \phi_1(x_n) & \ldots & \phi_n(x_n) \end{pmatrix} \begin{pmatrix} \overline{\phi}_1(x_1) & \ldots & \overline{\phi}_1(x_n) \\ \ldots & \ldots & \ldots \\ \overline{\phi}_n(x_1) & \ldots & \overline{\phi}_n(x_n) \end{pmatrix}. \quad (17)$$

In particular,

$$\det (K(x_i, x_j))_{1 \leq i,j \leq n} = \left( \sum_{\sigma \in S_n} \operatorname{sgn}(\sigma) \prod_{i=1}^{n} \phi_{\sigma_i}(x_i) \right) \left( \sum_{\tau \in S_n} \operatorname{sgn}(\tau) \prod_{i=1}^{n} \overline{\phi}_{\tau_i}(x_i) \right). \quad (18)$$

Now if $k_i$ are non-negative integers with $\sum_i k_i = n$, note that

$$\{\mathcal{X}(D_i) \geq k_i \text{ for all } 1 \leq i \leq r\} = \{\mathcal{X}(D_i) = k_i \text{ for all } 1 \leq i \leq r\},$$

since by Lemma 17, a determinantal process whose kernel defines a rank-$n$ projection operator has exactly $n$ points, almost surely. Thus, we have

$$\mathbf{P}\left[\mathcal{X}(D_i) = k_i \text{ for all } 1 \leq i \leq r\right] = \mathbf{E}\left[\prod_{i=1}^{r} \binom{\mathcal{X}(D_i)}{k_i}\right]$$
$$= \frac{1}{k_1! \cdots k_r!} \int_{\prod_{i=1}^{r} D_i^{k_i}} \det (K(x_k, x_\ell))_{1 \leq k, \ell \leq n} \, dx_1 \ldots dx_n$$
$$= \frac{1}{k_1! \cdots k_r!} \int_{\prod_{i=1}^{r} D_i^{k_i}} \sum_{\sigma, \tau} \operatorname{sgn}(\sigma) \operatorname{sgn}(\tau) \prod_{m=1}^{n} \phi_{\sigma_m}(x_m) \overline{\phi}_{\tau_m}(x_m) \, dx_1 \ldots dx_n.$$

Any term with $\sigma \neq \tau$ vanishes upon integrating. Indeed, if $\sigma(m) \neq \tau(m)$ for some $m$, then

$$\int_{D_{j(m)}} \phi_{\sigma_m}(x_m) \overline{\phi}_{\tau_m}(x_m) \, dx_m = 0$$

where $j(m)$ is the index for which

$$k_1 + \ldots + k_{j(m)-1} < m \leq k_1 + \ldots + k_{j(m)}.$$

Therefore,

$$\mathbf{E}\left[\prod_{i=1}^{r} \binom{\mathcal{X}(D_i)}{k_i}\right] = \frac{1}{k_1! \cdots k_r!} \sum_{\sigma} \prod_{m=1}^{n} \int_{D_{j(m)}} |\phi_{\sigma_m}(x)|^2 dx$$
$$= \frac{1}{k_1! \cdots k_r!} \sum_{\sigma} \prod_{m=1}^{n} \lambda_{\sigma_m, j(m)}.$$

Now consider (8) and set $M_i = \sum_k \xi_{k,i}$ for $1 \le i \le r$. Our goal is to compute $\mathbf{P}[M_1 = k_1, \ldots, M_r = k_r]$. This problem is the same as putting $n$ ball into $r$ cells, where the probability for the $j^{\text{th}}$ ball to fall in cell $i$ is $\lambda_{j,i}$. To have $k_i$ balls in cell $i$ for each $i$, we first take a permutation $\sigma$ of $\{1, 2, \ldots, n\}$ and then put the $\sigma_m^{\text{th}}$ ball into cell $j(m)$ if $k_1 + \ldots + k_{j(m)-1} < m \le k_1 + \ldots + k_{j(m)}$. However, this counts each assignment of balls $\prod_{i=1}^r k_i!$ times. This implies that

$$\mathbf{P}[M_1 = k_1, \ldots, M_r = k_r] = \frac{1}{k_1! \cdots k_r!} \sum_\sigma \prod_{m=1}^n \lambda_{\sigma_m, j(m)}.$$

Thus,
$$(\mathcal{X}(D_1), \ldots, \mathcal{X}(D_r)) \stackrel{d}{=} (M_1, \ldots, M_r), \tag{19}$$

which is precisely what we wanted to show.

Now we deal with the two assumptions that we made at the beginning. If $\cup_i D_i \ne \Lambda$, we could restrict the point process to $\cup_i D_i$. We still have a determinantal process. Then, if the kernel does not define a projection, apply Theorem 7 to write $\mathcal{X}$ as a mixture of determinantal projection processes. Applying (19) to each component in the mixture we obtain the theorem. The finite rank assumption can be relaxed in the same way as in Theorem 7. □

### *Applications*

As an application of Theorem 7 we can derive the following central limit theorem for determinantal processes due to Costin and Lebowitz [5] in case of the sine kernel, and due to Soshnikov [28] for general determinantal processes.

**Theorem 24.** *Let $\mathcal{X}_n$ be a sequence of determinantal processes on $\Lambda$ with kernels $K_n$. Let $D_n$ be a sequence of measurable subsets of $\Lambda$ such that $Var(\mathcal{X}_n(D_n)) \to \infty$ as $n \to \infty$. Then*

$$\frac{\mathcal{X}_n(D_n) - \mathbf{E}[\mathcal{X}_n(D_n)]}{\sqrt{Var(\mathcal{X}_n(D_n))}} \xrightarrow{d} N(0, 1). \tag{20}$$

*Proof.* By Theorem 7, $\mathcal{X}(D_n)$ has the same distribution as a sum of independent Bernoullis with parameters being the eigenvalues of the integral operators associated with $K_n$ restricted to $D$. A straightforward application of Lindeberg-Feller CLT for triangular arrays gives the result.

□

**Remark 25.** Existing proofs of results of the kind of Theorem 24 ([5], [28]) use the moment generating function for particle counts. Indeed, one standard way to prove central limit theorems (including the Lindeberg-Feller theorem) uses generating functions. The advantage of our proof is that the reason for the validity of the CLT is more transparent and a repetition of well known computations are avoided. Moreover, by applying the classical theory of sums of independent variables, local limit theorems, large deviation principles and extreme value asymptotics follow without any extra effort.

### *Radially symmetric processes on the complex plane*

Theorem 9 implies that when a determinantal process with kernel $K$ has the form $K(z,w) = \sum_n c_n(z\overline{w})^n$, with respect to a radially symmetric measure $\mu$, then the absolute values of the points are independent. More precisely, we have

**Theorem 26.** *Let $\mathcal{X}$ be a determinantal process with kernel $K$ with respect to a radially symmetric measure $\mu$ on $\mathbb{C}$. Write $K(z,w) = \sum_k \lambda_k a_k^2 (z\overline{w})^k$, where $a_k z^k$, $0 \leq k \leq n-1$ are the normalized eigenfunctions for $K$. The following construction describes the distribution of $\{|z|^2 : z \in \mathcal{X}\}$.*

- *Let $Z$ be picked from the probability distribution $\mu/\mu(\Lambda)$, and let $Q_0 = |Z|^2$.*
- *For $1 \leq k \leq n-1$ let $Q_k$ be an independent size-biased version of $Q_{k-1}$ (i.e., $Q_k$ has density $f_k(q) = \frac{a_k^2}{a_{k-1}^2} q$ with respect to the law of $Q_{k-1}$).*
- *Consider the point process in which each point $Q_k$ is included with probability $\lambda_k$ independently of everything else.*

When $\mu$ has density $\phi(|z|)$, then $Q_k$ has density

$$\pi a_k^2 q^k \phi(\sqrt{q}). \tag{21}$$

Theorem 26 (and its higher dimensional analogues) is the only kind of example that we know for interesting simultaneously observable counts.

*Proof.* Let $\nu$ be the measure of the squared modulus of a point picked from $\mu$. In particular, if $\mu$ has density $\phi(|z|)$, then we have $d\nu(q) = \pi\phi(\sqrt{q})\,dq$.

For $1 \leq i \leq r$, let $D_i$ be mutually disjoint open annuli centered at 0 with inner and outer radii $r_i$ and $R_i$ respectively. Since the functions $z^k$ are orthogonal on any annulus centered at zero, it follows that the $D_i$'s are simultaneously observable. To compute the eigenvalues, we integrate these functions against the restricted kernel; clearly, all terms but one cancel, and we get that for $z \in D_i$

$$z^k \lambda_{k,i} = \int_{D_i} \lambda_k a_k^2 (z\overline{w})^k w^k d\mu(w), \quad \text{and so}$$

$$\lambda_{k,i} = \lambda_k a_k^2 \int_{D_i} |w|^{2k} d\mu(w)$$

$$= \lambda_k a_k^2 \int_{r_i^2}^{R_i^2} q^k d\nu(q).$$

As $r_i, R_i$ change, the last expression remains proportional to the probability that the $k$ times size-biased random variable $Q_k$ falls in $(r_i^2, R_i^2)$. When we set $(r_i, R_i) = (0, \infty)$, the result is $\lambda_k$ because $a_k w^k$ has norm 1. Thus the constant of proportionality equals $\lambda_k$. The theorem now follows from Proposition 9. □

**Example 27 (Ginibre Ensemble revisited).** Recall that the $n^{\text{th}}$ Ginibre ensemble described in Example 15 is the determinantal process $\mathbf{G}_n$ on $\mathbb{C}$ with kernel $K_n(z,w) = \sum_{k=0}^{n-1} \lambda_k a_k^2 (z\overline{w})^k$ with respect to the complex Gaussian measure $d\mu = \frac{1}{\pi} e^{-|z|^2} dz$, where $a_k^2 = 1/k!$, and $\lambda_k = 1$. The modulus-squared of a complex Gaussian is a gamma$(1, 1)$ random variable, and its $k$-times size-biased version has gamma$(k + 1, 1)$ distribution (see (21)). Theorem 26 immediately yields the following.

**Theorem 28 (Kostlan [17]).** *The set of absolute values of the points of $\mathbf{G}_n$ has the same distribution as $\{Y_1, \ldots, Y_n\}$ where $Y_i$ are independent and $Y_i^2 \sim$ gamma$(i, 1)$.*

All of the above holds for $n = \infty$ also, in which case we have a determinantal process with kernel $e^{z\overline{w}}$ with respect to $d\mu = \frac{1}{\pi} e^{-\frac{1}{2}|z|^2} dz$. This case is also of interest as $\mathbf{G}_\infty$ is a translation invariant process in the plane.

**Example 29 (Zero set of a Gaussian analytic function).** Recall the zero set of $\mathbf{f}_1(z) := \sum_{n=0}^{\infty} a_n z^n$ is a determinantal process in the disk with the Bergman kernel

$$K(z, w) = \frac{1}{\pi(1 - z\overline{w})^2} = \frac{1}{\pi} \sum_{k=0}^{\infty} (k + 1)(z\overline{w})^k,$$

with respect to Lebesgue measure in the unit disk. Theorem 26 applies with $a_k^2 = (k + 1)$ and $\lambda_k = 1$ (to make $K$ trace class, we first have to restrict it to the disk of radius $r < 1$ and let $r \to 1$). From (21) we immediately see that $Q_k$ has beta$(k + 1, 1)$ distribution. Equivalently, we get the following.

**Theorem 30 (Peres and Virág [23]).** *The set of absolute values of the points in the zero set of $\mathbf{f}_1$ has the same distribution as $\{U_1^{1/2}, U_2^{1/4}, U_3^{1/6}, \ldots\}$ where $U_i$ are i.i.d. uniform$[0, 1]$ random variables.*

We can of course consider the determinantal process with the truncated Bergman kernel $K_n(z, w) = \frac{1}{\pi} \sum_{k=0}^{n-1} (k + 1)(z\overline{w})^k$. The set of absolute values of this process has the same distribution as $\{U_1^{1/2}, \ldots, U_n^{1/2n}\}$.

## 3. High powers of complex polynomial processes

Rains [24] showed that sufficiently high powers of eigenvalues of a random unitary matrix are independent.

**Theorem 31 (Rains [24]).** *Let $\{z_1, \ldots, z_n\}$ be the set of eigenvalues of a random unitary matrix chosen according to Haar measure on $\mathcal{U}(n)$. Then for every $k \geq n$, $\{z_1^k, \ldots, z_n^k\}$ has the same distribution as a set of $n$ points chosen independently according to uniform measure on the unit circle in the complex plane.*

We point out that this theorem holds whenever the angular distribution of the points is a trigonometric polynomial.

**Proposition 32.** *Let $\{z_1, \ldots, z_n\}$ be distributed on $(S^1)^{\otimes n}$ with density $P(z_1, \ldots, z_n, \overline{z}_1, \ldots, \overline{z}_n)$ w.r.t. uniform measure on $(S^1)^{\otimes n}$, where $P$ is a polynomial of degree $d$ or less in each variable. Then for every $k > d$ the vector* $(z_1^k, \ldots, z_n^k)$ *has the distribution of $n$ points chosen independently according to uniform measure on $S^1$.*

*Proof.* Fix $k > d$ and consider any joint moment of $(z_1^k, \ldots, z_n^k)$,

$$\mathbf{E}\left[\prod_{i=1}^n \left(z_i^{km_i} \overline{z}_i^{k\ell_i}\right)\right] = \int_{(S^1)^{\otimes n}} \prod_{i=1}^n \left(z_i^{km_i} \overline{z}_i^{k\ell_i}\right) P(z_1, \ldots, z_n, \overline{z}_1, \ldots, \overline{z}_n) d\lambda,$$

where $\lambda$ denotes the uniform measure on $(S^1)^{\otimes n}$. If $m_i \neq \ell_i$ for some $i$ then the integral vanishes. To see this, note that the average of a monomial over $(S^1)^{\otimes n}$ is either 1 or 0 depending on whether the exponent of every $z_i$ matches that of $\overline{z}_i$. Suppose without loss of generality that $m_1 > \ell_1$. Then in each term, we have an excess of $z_1^k$ which cannot be matched by an equal power of $\overline{z}_1$ because $P$ has degree less than $k$ as a polynomial in $\overline{z}_1$.

We conclude that the joint moments are zero unless $m_i = \ell_i$ for all $i$. If $m_i = \ell_i$ for all $i$, then the expectation equals 1. Thus, the joint moments of $(z_1^k, \ldots, z_n^k)$ are the same as those of $n$ i.i.d. points chosen uniformly on the unit circle. This proves the proposition. □

More generally, by conditioning on the absolute values we get the following.

**Corollary 33.** *Let $\zeta_1, \ldots, \zeta_k$ be complex random variables with distribution given by*

$$P(z_1, \ldots, z_n, \overline{z}_1, \ldots, \overline{z}_n) d\mu_1(z_1) \cdots d\mu_n(z_n),$$

*where $P$ is a polynomial of degree $d$ or less in each variable, and the measures $\mu_i$ are radially symmetric. Then for every $k > d$, the angles $\mathrm{Arg}(\zeta_1^k), \ldots, \mathrm{Arg}(\zeta_n^k)$ are independent, have uniform distribution, and are independent of the moduli $\{|\zeta_1|, \ldots, |\zeta_n|\}$.*

Corollary 33 applies to powers of points of determinantal processes with kernels of the form $K(z, w) = \sum_{k=0}^d c_k (z\overline{w})^k$ w.r.t a radially symmetric measure $\mu$ on the complex plane. Combining this observation with our earlier results on the independence of the absolute values of the points, we get the following result.

**Theorem 34.** *Let $\mathcal{X} = \{z_1, \ldots, z_n\}$ be a determinantal process on the complex plane with kernel $K(z, w) = \sum_{\ell=0}^d c_k (z\overline{w})^\ell$ (satisfying $\#\{k \geq 0 : c_k \neq 0\} = n$) with respect to a radially symmetric measure $\mu$. Then for every $\ell \geq d$, the points $\{z_1^\ell, \ldots, z_n^\ell\}$ are distributed as a set of independent random variables with radially symmetric distribution.*

## 4. Permanental processes

In this section we prove analogous theorems for permanental processes. We begin with the following known representation of permanental processes, that can be found in Macchi [22].

**Proposition 35.** *Let $F$ be a complex Gaussian process on $\Lambda$. Given $F$, let $\mathcal{X}$ be a Poisson process in $\Lambda$ with intensity $|F|^2$. Then $\mathcal{X}$ is a permanental process with kernel $K(x, y) = \mathbf{E}\left[F(x)\overline{F}(y)\right]$.*

**Remark 36.** Since any non-negative definite Hermitian kernel is the covariance kernel of a complex Gaussian process, it follows that all permanental processes are of the above form.

*Proof.* Given $F$, the joint intensities of $\mathcal{X}$ are $\tilde{\rho}_k(x_1, \ldots, x_k) = \prod_{i=1}^{k} |F(x_i)|^2$. Hence it follows that the unconditional joint intensities of $\mathcal{X}$ are $\rho_k(x_1, \ldots, x_k) = \mathbf{E}\left[\prod_{i=1}^{k} |F(x_i)|^2\right]$. Now apply the classical Wick formula (see [11] or [23], §2, Fact 10) to get the result. □

**Corollary 37.** *If $K$ determines a self-adjoint non-negative definite locally trace-class integral operator $\mathcal{K}$, then there exists a permanental process with kernel $K$.*

Now we prove Theorem 10 using the representation in Proposition 35. We need the following simple fact:

**Fact 38.** Let $\mathcal{Y}$ be a Poisson process on $\Lambda$ with intensity measure $\nu$. Assume that $\nu(\Lambda) < \infty$ and $\nu$ is absolutely continuous with respect to $\mu$. Let $Y$ be the random vector of length $\mathcal{Y}(\Lambda)$ obtained from a uniform random ordering of the points of $\mathcal{Y}$. For $k \geq 1$, the law of $Y$ on the event that $\mathcal{Y}(\Lambda) = k$ is a subprobability measure on $\Lambda^k$ with density

$$g_k(z_1, \ldots, z_k) = \frac{1}{k!} \left[ e^{-\nu(\Lambda)} \prod_{i=1}^{k} \frac{d\nu}{d\mu}(z_i) \right] \quad (22)$$

with respect to $\mu^k$. We have $\int g_k \, d\mu^k = \mathbf{P}\left[\mathcal{Y}(\Lambda) = k\right]$.

*Proof of Theorem 10:.* We use the construction in Proposition 35 with $F(z) = \sum_{k=1}^{n} \sqrt{\lambda_k} a_k \phi_k(z)$ where $a_k$ are independent standard complex Gaussian random variables. Let $X$ be the random vector obtained from a uniform random ordering of the points of $\mathcal{X}$. If we first condition on $F$, then by Fact 38 the joint density of the random vector $X$ on the event $\{\mathcal{X}(\Lambda) = k\}$ is given by

$$j_{F,k}(z_1, \ldots, z_k) = \frac{1}{k!} \left[ e^{-\int_\Lambda |F(x)|^2 d\mu(x)} \prod_{i=1}^{k} |F(z_i)|^2 \right],$$

which is a subprobability measure with total weight $\mathbf{P}\left[\mathcal{X}(\Lambda) = k \,\big|\, F\right]$. Integrating over the distribution of $F$ we get that on the event $\{\mathcal{X}(\Lambda) = k\}$ the density of $X$ is

$$\begin{aligned} j_k(z_1, \ldots, z_k) &= \frac{1}{k!} \mathbf{E}\left[ e^{-\int_\Lambda |F(x)|^2 d\mu(x)} \prod_{i=1}^{k} |F(z_i)|^2 \right] \\ &= \frac{1}{k!} \mathbf{E}\left[ e^{-\sum_m \lambda_m |a_m|^2} \left| \prod_{i=1}^{k} \left( \sum_m \sqrt{\lambda_m} a_m \phi_m(z_i) \right) \right|^2 \right], \quad (23) \end{aligned}$$

which is also a subprobability measure with total weight $\mathbf{P}\left[\mathcal{X}(\Lambda) = k\right]$. We now expand the product inside the expectation (23) as a sum indexed by ordered set partitions $(S_1, \ldots, S_n)$ and $(T_1, \ldots, T_n)$ of $\{1, 2, \ldots, k\}$. The set partitions corresponding to a summand $q$ are constructed by letting $S_\ell$ be the set of indices $i$ for which $q$ contains the term $\sqrt{\lambda_\ell} a_\ell \phi_\ell(z_i)$ and $T_i$ be the set of indices $i$ for which $q$ contains the term $\sqrt{\lambda_\ell}\overline{a}_\ell\overline{\phi}_\ell(z_i)$. The summand corresponding to the partitions $(S_\ell), (T_\ell)$ is thus:

$$\mathbf{E}\left[e^{-\sum_m \lambda_m|a_m|^2}\left(\prod_\ell a_\ell^{|S_\ell|}\overline{a}_\ell^{|T_\ell|}\right)\left(\prod_\ell \prod_{i\in S_\ell}\lambda_\ell^{|S_\ell|/2}\phi_\ell(z_i)\right)\right.$$
$$\left.\times\left(\prod_\ell\prod_{i\in T_\ell}\lambda_\ell^{|T_\ell|/2}\overline{\phi}_\ell(z_i)\right)\right],$$

which clearly vanishes unless $|S_\ell| = |T_\ell|$ for every $\ell$. Also note that for a standard complex normal random variable $a$,

$$\mathbf{E}\left[e^{-\lambda|a|^2}|a|^{2m}\right] = \frac{m!}{(1+\lambda)^{m+1}}.$$

Therefore by fixing an ordered partition of the integer $k$ with $n$ parts (some of the parts may be zero) and then summing over all ordered set partitions $(S_\ell), (T_\ell)$ with those sizes, we find that

$$j_k(z_1,\ldots,z_k) = \frac{1}{k!}\sum_{(m_1,\ldots,m_n):\sum m_i=k}\prod_{i=1}^n \frac{\lambda_i^{m_i}m_i!}{(1+\lambda_i)^{m_i+1}}$$
$$\times\left|\sum_{(S_i):|S_i|=m_i}\prod_{\ell\geq 1}\prod_{i\in S_\ell}\phi_\ell(z_i)\right|^2. \qquad (24)$$

Now it is easy to see that

$$\sum_{(S_i):|S_i|=m_i}\prod_{\ell\geq 1}\prod_{i\in S_\ell}\phi_\ell(z_i)$$
$$= \left(\prod_{i=1}^n \frac{1}{m_i!}\right)\mathrm{per}\begin{bmatrix}\{\phi_1(z_1) & \cdots & \phi_1(z_k)\}\,m_1 \\ \cdots & \cdots & \cdots \\ \{\phi_n(z_1) & \cdots & \phi_n(z_k)\}\,m_n\end{bmatrix}.$$

Recall that the geometric random variables $\chi_i$ in the statement of the theorem have the distributions $\mathbf{P}[\chi_i = m] = \frac{\lambda_i^m}{(1+\lambda_i)^{1+m}}$. Therefore we obtain

$$j_k(z_1,\ldots,z_k) = \sum_{(m_1,\ldots,m_n):\sum m_i=k}\prod_{i=1}^n \frac{\mathbf{P}[\chi_i=m_i]}{k!\prod_{i=1}^n m_i!}$$
$$\times\left|\mathrm{per}\begin{bmatrix}\{\phi_1(z_1) & \cdots & \phi_1(z_k)\}\,m_1 \\ \cdots & \cdots & \cdots \\ \{\phi_n(z_1) & \cdots & \phi_n(z_k)\}\,m_n\end{bmatrix}\right|^2. \qquad (25)$$

Now we integrate (25) over all the variables $z_i$. Write the absolute square of the permanent on the right as

$$\mathrm{per}\begin{bmatrix}\{\phi_1(z_1) & \cdots & \phi_1(z_k)\}\,m_1 \\ \cdots & \cdots & \cdots \\ \{\phi_n(z_1) & \cdots & \phi_n(z_k)\}\,m_n\end{bmatrix}\overline{\mathrm{per}}\begin{bmatrix}\{\phi_1(z_1) & \cdots & \phi_1(z_k)\}\,m_1 \\ \cdots & \cdots & \cdots \\ \{\phi_n(z_1) & \cdots & \phi_n(z_k)\}\,m_n\end{bmatrix}$$

and expand these two factors over permuatations $\pi$ and $\sigma$ of $\{1, 2, \ldots, k\}$. Letting $I_j$ denote the interval of integers $\{1 + \sum_{r=1}^{j-1} m_r, \ldots, \sum_{r=1}^{j} m_r\}$ we get a sum of terms of the form

$$\left(\prod_{j=1}^{n} \prod_{i \in \pi^{-1}(I_j)} \phi_j(z_i)\right) \left(\prod_{j=1}^{n} \prod_{i \in \sigma^{-1}(I_j)} \overline{\phi}_j(z_i)\right).$$

By orthogonality of $\phi_j$s, this term vanishes upon integration unless $\pi^{-1}(I_j) = \sigma^{-1}(I_j)$ for every $1 \leq j \leq n$. For a given $\pi$, there are $\prod_{j=1}^{n} m_j!$ choices of $\sigma$ that satisfy this. For each such $\sigma$, we get 1 upon integration over $z_i$s. Summing over all $k!$ choices for $\pi$, we get

$$\int_{\Lambda^k} j_k \, d\mu = \mathbf{P}\left[\mathcal{X}(\Lambda) = k\right] = \sum_{(m_1, \ldots, m_n) : \sum m_i = k} \prod_{i=1}^{n} \mathbf{P}\left[\chi_i = m_i\right]$$

$$= \mathbf{P}\left[\sum_{i=1}^{n} \chi_i = k\right],$$

which proves the claim about the number of points in $\Lambda$. Thus by (25) $\mathcal{X}$ is a mixture of the processes $\mathcal{X}_{\vec{m}}(D)$, with weights given by $\prod_{i=1}^{n} \mathbf{P}\left[\chi_i = m_i\right]$, where $\vec{m} = (m_1, \ldots, m_n)$ with $m_i$ being non-negative integers. This is what we wanted to prove. □

Now we prove Theorem 12. As before, we remark that it is applicable to the restriction of $\mathcal{X}$ to $D$ for any Borel set $D \subset \Lambda$.

*Proof of Theorem 12.* Suppose $D_1, \ldots, D_r$ are simultaneously observable as in the statement of the theorem. Use Proposition 35 to write $\mathcal{X}$ as a Poisson process with intensity $|F(x)|^2$ where $F$ is a Gaussian process with covariance kernel $K$. Explicitly, $F(x) = \sum_k a_k \sqrt{\lambda_k} \phi_k(x)$ for $x \in \cup_{i=1}^{r} D_i$, where $a_k$ are i.i.d. standard complex Gaussians, i.e., the real and imaginary parts of $a_k$ are i.i.d. $N(0, \frac{1}{2})$. Then given $\{a_k\}$, we know that $\mathcal{X}(D_i)$, $1 \leq i \leq r$ are independent Poisson($\int_{D_i} |F(x)|^2 d\mu(x)$). Writing $\int_{D_i} |F(x)|^2 d\mu(x) = \sum_k \lambda_{k,i} |a_k|^2$, we see that conditionally given $\{a_k\}$, the variables $\mathcal{X}(D_i)$ for $1 \leq i \leq r$ have the same distribution as $\sum_k (\eta_{k,1}, \ldots, \eta_{k,r})$, where $\{\eta_{k,i}\}_{1 \leq i \leq r}$ are chosen by sampling $\eta_k$ according to Poisson($\lambda_k |a_k|^2$) distribution and then assigning $\eta_k$ to the cells $D_i$ multinomially with probabilities $\frac{\lambda_{k,i}}{\lambda_k}$.

Now when we integrate over the randomness in $\{a_k\}$, we see that $\eta_k$ has a Geometric distribution with mean $\lambda_k$ and given $\eta_k$, the vector $(\eta_{k,1}, \ldots, \eta_{k,r})$ is still Multinomial($\eta_k; \frac{\lambda_{k,1}}{\lambda_k}, \ldots, \frac{\lambda_{k,r}}{\lambda_k}$). This completes the proof. □

### *Generalization: α-determinantal processes*

One way to generalize the concept of determinantal and permanental processes is to consider point processes with joint intensities given by

$$\rho_n(x_1, \ldots, x_n) = \det{}_\alpha(K(x_i, x_j)) \stackrel{def}{=} \sum_{\pi \in S_n} \alpha^{n - \nu(\pi)} \prod_{i=1}^{n} K(x_i, x_{\pi(i)}), \qquad (26)$$

where $\nu(\pi)$ is the number of cycles in the permutation $\pi$.

Such point processes are called $\alpha$-determinantal processes and were introduced by Vere-Jones in [29]. The values $\alpha = -1$ and $\alpha = +1$ correspond to determinantal and permanental processes, respectively. It is easy to check that the proof of Theorem 7 can be modified to get:

**Proposition 39.** *For a point process with joint intensities given by (26) (when it exists), the number of points that fall in any subregion is:*

- *a sum of independent Binomial$(-\frac{1}{\alpha}, -\alpha\lambda_k)$ random variables, if $-\frac{1}{\alpha}$ is a positive integer.*
- *a sum of independent Negative Binomial$(\frac{1}{\alpha}, \frac{\alpha\lambda_k}{\alpha\lambda_k+1})$ random variables, if $\alpha > 0$.*

In fact, if $-\frac{1}{\alpha}$ is a positive integer, this process is just a union of $-\frac{1}{\alpha}$ i.i.d. copies of the determinantal process with kernel $-\alpha K$. Similarly, if $\frac{1}{\alpha}$ is a positive integer, this process is a union of $\frac{1}{\alpha}$ i.i.d. copies of the permanental process with kernel $\alpha K$. More generally, the union of $m$ i.i.d. copies of the process corresponding to $\alpha$ and kernel $K$ gives a process distributed according to $\frac{\alpha}{m}$ and kernel $mK$. If $K$ is real, then $\frac{2}{m}$-determinantal processes also exist [25]. For $\alpha > 0$, little is known about the existence of $\alpha$-determinantal processes beyond these examples. Shirai and Takahashi [25] conjecture the following:

**Conjecture 40.** *If $K$ is a (real) kernel defining a self-adjoint, non-negative integral operator on $L^2(\Lambda)$ and $0 \leq \alpha \leq 2$, then the $\alpha$-determinantal process with kernel $K$ exists. However, if $\alpha > 2$, then there exists such a kernel $K$ for which there is no corresponding $\alpha$-determinantal process.*

We verify this conjecture for $\alpha > 4$. Let $\Lambda$ be a discrete space consisting of three points, and consider the $3 \times 3$ matrix

$$K = \begin{pmatrix} 2 & -1 & -1 \\ -1 & 2 & -1 \\ -1 & -1 & 2 \end{pmatrix}.$$

It is easy to check that the eigenvalues of $K$ are $3, 3, 0$ and

$$\det\nolimits_\alpha(K(i,j))_{1 \leq i,j \leq 3} = 2(4-\alpha)(\alpha+1),$$

which is negative for $\alpha > 4$. Since point processes must have non-negative joint intensities, we conclude that no $\alpha$-determinantal processes with this kernel can exist for $\alpha > 4$.

**Remark 41.** It was pointed out to us by Steve Evans that this question is closely related to the question of infinite divisibility of a multivariate version of Gamma distribution that has been studied by Griffiths [12] and Griffiths and Milne [13]. In short, knowing that a certain multivariate Gamma distribution is *not* infinitely divisible would tell us that $\alpha$-permanental processes for a related kernel do *not* exist for sufficiently large $\alpha$ (but does not give us the smallest $\alpha$ for which this happens).

**Acknowledgements.** The insightful recent papers by Soshnikov [27], Shirai-Takahashi [25] and Lyons [20] were of great help in preparing this work. We are grateful to Russell Lyons and Balázs Szegedy for useful conversations, and to Persi Diaconis, Alice Guionnet and Ofer Zeitouni for comments. We thank Steve Evans for several comments and for pointing out the references to $\alpha$-permanental processes and papers on multivariate Gamma distributions.